\documentclass[10pt]{amsart}

\usepackage{amssymb,eucal}
\usepackage{amscd}
\usepackage[all,cmtip]{xy}
\usepackage{rotating}
\usepackage{hyperref,mathrsfs,soul,upgreek,mdframed}

\setul{}{1.3pt}

\usepackage{graphicx,xcolor}
\usepackage{latexsym}
\usepackage{amssymb,amsmath,array,delarray,amsfonts,times}

\newcommand{\eop}{\hspace*{\fill}$\blacksquare$}
\newcommand{\I}{\mathbf{I}}

\newcommand{\btt}{\begin{tTheorem}}
\newcommand{\ett}{\end{tTheorem}}

\newcommand{\Aut}{\mathrm{Aut}}
\newcommand{\bt}{\begin{Theorem}}
\newcommand{\et}{\end{Theorem}}
\newcommand{\bmt}{\begin{mainTheorem}}
\newcommand{\emt}{\end{mainTheorem}}
\newcommand{\bc}{\begin{Corollary}}
\newcommand{\bl}{\begin{Lemma}}
\newcommand{\ec}{\end{Corollary}}
\newcommand{\el}{\end{Lemma}}
\newcommand{\bo}{\begin{observation}}
\newcommand{\eo}{\end{observation}}
\newcommand{\bp}{\begin{proposition}}
\newcommand{\ep}{\end{proposition}}
\newcommand{\br}{\begin{remark}}
\newcommand{\er}{\end{remark}}

\newcommand{\bcc}{\begin{conjecture}}
\newcommand{\ecc}{\end{conjecture}}

\newcommand{\mS}{\ensuremath{\mathcal{S}}}
\newcommand{\mO}{\ensuremath{\mathcal{O}}}

\newcommand{\mL}{\ensuremath{\mathcal{L}}}
\newcommand{\mT}{\ensuremath{\mathcal{T}}}

\newcommand{\mC}{\ensuremath{\mathcal{C}}}
\newcommand{\mQ}{\ensuremath{\mathcal{Q}}}

\newcommand{\mP}{\ensuremath{\mathcal{P}}}
\newcommand{\mB}{\ensuremath{\mathcal{B}}}
\newcommand{\mX}{\ensuremath{\mathcal{X}}}
\newcommand{\mI}{\mathcal{I}}

\newtheorem{Theorem}{Theorem}[section]
\newtheorem{observation}[Theorem]{Observation}
\newtheorem{Corollary}[Theorem]{Corollary}

\newtheorem{Lemma}[Theorem]{Lemma}
\newtheorem{proposition}[Theorem]{Proposition}
\newtheorem{remark}[Theorem]{Remark}

\newtheorem{conjecture}[Theorem]{Conjecture}

\newcommand{\bT}{\texttt{T}}
\newcommand{\bI}{\texttt{I}}

\newpage

\title[Two questions of Moorhouse]{Two questions of Moorhouse on indiscernible locally finite\\ generalized quadrangles}

\author[K. Thas]{Koen Thas}

\address{{Ghent University},
{Department of Mathematics: algebra \& geometry},
{Krijgslaan 281, S25, B-9000 Ghent, Belgium}}

\date{May 2023}

\email{koen.thas@gmail.com}

\urladdr{http://cage.ugent.be/$\sim$kthas}

\begin{document}
\maketitle

\begin{abstract}
We settle two questions posed by G. Eric Moorhouse on the existence of locally finite generalized quadrangles with indiscernible ovoids or spreads.  
\end{abstract}

\tableofcontents

\bigskip
\section{Introduction}

Consider a thick generalized $n$-gon $\Upgamma$ \cite{HVM} with $s+1$ points on each line and $t+1$ lines through each point. If $n$ is odd, then it is easy to show that $s=t$, see \cite[section 1.5.3]{HVM}. If $n$ is even though, there are examples where $s\neq t$, a most striking example being $n=8$ in which case a theorem of Feit and Higman \cite{FeHi} implies that if $st$ is finite, $2st$ is a perfect square and so $s$ is never equal to $t$. If both $s$ and $t$ are finite, they are bounded by each other; to be more specific, $s\leq t^2\leq s^4$ for $n=4$ and $n=8$ by results of Higman \cite{Hi75} (1975),
and $s\leq t^3\leq s^9$ for $n = 6$ by Haemers and Roos \cite{HR} (1981). These results can also be found in \cite[section 1.7.2]{HVM}.
For other even values of $n$, $\Upgamma$ cannot exist by a famous result of Feit and Higman which also appeared  in their 1964 paper \cite{FeHi}. In {\em loc. cit.}, 
several necessary divisibility conditions involving the parameters $s$ and $t$ of a generalized $n$-gon can be found (with $n \in \{ 3,4,6,8\}$), in order for a generalized $n$-gon with these parameters to exist. \\

The very first question which arises about parameters of generalized polygons | before one is ready to formulate prime power conjectures in the finite case | obviously is: ``Do thick generalized polygons exist with \ul{one finite parameter} and \ul{one non-finite parameter}?'' First posed by Tits in the 1960s, this question remains a mystery in our knowledge about the fundaments of generalized polygons. \\

By definition, we call such polygons {\em locally finite}.  Note that in Van Maldeghem's book \cite{HVM}, such generalized polygons are called {\em semi-finite}.\\
%Tits's question can be found as Problem 5 in the ``Ten Most Famous Open Problems'' chapter of 
%Van Maldeghem's book \cite{HVM}, see also \cite[section 10]{JATGP}.)

Not much is known on Tits's question. All the known results comprise the case $n = 4$. Here they are.

\begin{itemize}
\item
P. J. Cameron \cite{PJC} showed in 1981 that if $n=4$ and $s=2$, then $t$ is finite. 
\item
In \cite{AEB} A. E. Brouwer shows the same thing for $n=4$ and $s=3$, and the proof is purely combinatorial (unlike a nonpublished but earlier proof of Kantor \cite{Kanunp}). 
\item
More recently, G. Cherlin used Model Theory (in \cite{Cher}) to handle the generalized $4$-gons  with five points on a line.
\item
For other values of $n$ and $s$ (where $n$ of course is even), nothing is known without any extra assumptions. 
\end{itemize}

Apart from the aforementioned results, there is only one other ``general result'' on parameters of generalized polygons (without invoking additional structure through, e.g., the existence of certain substructures or the occurrence of certain group actions):

\bt[Bruck and Ryser \cite{BR}, 1949]
If $\Upgamma$ is a finite projective plane of order $m$, $m \in \mathbb{N}^\times$, and $m \equiv 1,2\mod{4}$, then $m$ is the sum of two perfect squares.
\et

Note that it is not hard to construct generalized quadrangles (= generalized $4$-gons) with parameters $(\upalpha,\upbeta)$, where $\upalpha$ and $\upbeta$ are different cardinal numbers both not a finite positive integer.
Let $I$ be an infinite uncountable set, and consider the vector space $\mathbb{Q}^{\vert I\vert}$ consisting of all $\vert I\vert$-tuples $(q_i)_{i \in I}$, $q_i \in \mathbb{Q}$, with only a 
finite number of nonzero entries. We assume without loss of generality that $I$ contains the symbols $0$ and $1$.
Define the following quadratic form:
\begin{equation}
\iota: X_0^2 + X_1^2 - \sum_{i \in I \setminus \{0,1\}}X_i^2.
\end{equation}

Then $\iota$ has Witt index $2$, and
the corresponding classical orthogonal quadrangle $\mQ(\vert I\vert,\mathbb{Q},\iota)$ (cf. \cite[chapter 2]{HVM})
is a Moufang generalized quadrangle with $\vert \mathbb{Q}\vert$ points per line and $\vert I\vert$ lines on a point. It is fully embedded in the projective space $\mathbf{P}(\mathbb{Q}^{\vert I\vert})$ (cf. \cite[chapter 2]{HVM}). Other examples of buildings with ``mixed parameters'' can be found in \cite{Order}.

\medskip
\subsection*{This paper}

This note is inspired by G. Eric Moorhouse's notes \cite{Moor} on a paper of Cherlin \cite{Cher} which deals with locally finite quadrangles with $5$ points per line. Moorhouse's text ends 
with two basic questions (stated in section \ref{Moorques} of this paper) on the model theory of hypothetical locally finite generalized quadrangles, and answering these questions is the main purpose of this paper.

In \cite{Cher}, Cherlin observes that if a locally finite generalized quadrangle would happen to exist, then a generalized quadrangle with the same parameters would exist with an infinite and independent set of points or lines (that is, 
an infinite partial ovoid or partial spread) which are in a precise sense coordinatized over some totally ordered set $(S,\leq)$, and
with extreme symmetry properties relative to the automorphism group of the quadrangle. Such sets are called ``indiscernible'' over $(S,\leq)$. 
It should be emphasized that one is allowed to choose $(S,\leq)$ arbitrary in this result. 
Moorhouse's questions amount to considering the extremal cases 
for these sets: {\em ovoids} and {\em spreads}. 

In this note, we will show that locally finite generalized quadrangles with ovoids or spreads which are indiscernible over any given infinite totally ordered set $(S,\leq)$, cannot exist. As we will indicate below, different types of indiscernibility will be considered, depending on which case (ovoid/spread) we are handling. 
These results should be seen as a very first step in attacking Tits's problem through Model Theory, after Cherlin's seminal note \cite{Cher}. 
A first version of the proof benefited greatly from the proof of an old result by Dushnik and Miller \cite{DUMI} on similarity transformations of totally ordered sets, which dates back to 1940. When the paper was almost completed, the author noticed a geometrically much  simpler proof which only requires the existence of certain isomorphisms  that  appear to be naturally present in abundance under Moorhouse's conditions. Also, initially only quadrangles were handled with a countable number of points and lines. This restriction is no longer present. 

%As we will see, the structure of such sets has great structural implications for generalized quadrangles. 

%\br
%{\rm
%Dozens of surveys have been written on parameter problems for generalized polygons, and some of these problems are considered to be amongst the most fundamental in the entire theory, even when %invoking heavy extra assumptions (such as a sharply transitive automorphism  group action on the points or flags, the occurrence of a group with a BN-pair, the occurrence of a translation group, etc.). And %one also has the famous {\em Prime Power Conjecture} for projective planes, which states that the order of a finite projective plane is always a prime power. We refer the reader to the recent survey %\cite{Order} of the author, which 
%investigates ``prime power conjectures'' in a very general context, and also to the many references therein.\\  
%}
%\er

%\br[Burnside Polygons]
%{\rm The Burnside problem, posed by William Burnside in 1902 and one of the oldest and most influential questions in Group Theory,
%asks whether a finitely generated group in which every element has finite order must necessarily be a finite group. Such infinite groups indeed exist, and are usually called {\em Burnside groups}. Since the %generalized polygons we encounter are geometric analogues of such groups, the term ``Burnside polygon'' seems in place. For more on Burnside's problem, we refer to \cite{Zelm}. Note that 
%since Burnside groups do exist, a straightforward adaptation of Theorem \ref{Cherlem} below holds true for this class of groups.}\\
%\er 

%(Zelmanov. BURNSIDE HVMYGONS.)

\subsection*{ACKNOWLEDGMENT}

I would like to thank Gregory Cherlin, Alex Kruckman and Eric Moorhouse for several valuable communications on the topic of this paper.  \\

\bigskip
\section{Some notes on generalized quadrangles}

In this paper, we solely concentrate on {\em generalized quadrangles}. \\

A {\em generalized $4$-gon}, also {\em generalized quadrangle} (abbreviated as ``GQ''), is a point-line incidence geometry $\Upgamma = (\mP,\mB,\I)$ for which  the following axioms are satisfied:

\begin{itemize}
\item[(i)]
$\Upgamma$ contains no ordinary $k$-gon (as a subgeometry), for $2 \leq k < 4$;
\item[(ii)]
any two elements $x,y \in \mP \cup \mB$ are contained in some ordinary $4$-gon in
$\Upgamma$;
\item[(iii)]
there exists an ordinary $5$-gon in $\Upgamma$.\\
\end{itemize}

The value $4$ is called the {\em gonality} of the GQ. 

By (iii), generalized quadrangles have at least three points per line and three lines per point. So by this definition, we do not consider ``thin'' quadrangles.  Note that points and lines play the same role in the defining axioms; this is the principle of ``duality''.\\

It can be shown that generalized quadrangles have an {\em order} $(s,t)$; there exist constants $s, t$ such that 
the number of points incident with a line is $s + 1$, and the number of lines incident  with a point is 
$t + 1$, cf. \cite[section 1.5.3]{HVM}. \\

A {\em sub generalized quadrangle} or {\em subquadrangle} of a GQ $\Upgamma = (\mP,\mB,\I)$ is a GQ $\Upgamma' = (\mP',\mB',\I')$ for which $\mP \subseteq \mP$,
$\mB' \subseteq \mB$ and $\I' \subseteq \I$.  It is {\em full} if for any line $L$ of $\Upgamma'$, we have that $x \I'  L$ if and only if $x \I L$. Dually, we define {\em ideal subquadrangles}.\\

{ 
Suppose $x$ is a point in a generalized quadrangle $\Upgamma$. Then $x^{\perp}$ is by definition the set of points in $\Upgamma$ which are collinear with $x$. Now if $A$ is more generally a point set of $\Upgamma$, then $A^{\perp} = \bigcap_{a \in A}a^\perp$, and $A^{\perp\perp} = {\Big( A^\perp \Big)}^{\perp}$. 
Dually, we use the same notation for lines.}\\

An {\em automorphism} of a generalized quadrangle $\Upgamma = (\mP,\mB,\I)$ is a bijection of $\mP \cup \mB$ which preserves $\mP$, $\mB$ and incidence.
The full set of automorphisms of a GQ forms a group in a natural way | the {\em automorphism group} of $\Upgamma$, denoted $\Aut(\Upgamma)$. It is one of its most important invariants.
If $B$ is an automorphism group of a generalized quadrangle $\Upgamma = (\mP,\mB,\I)$, and $R$ is a subset of $\mP$, $B_{[R]}$ is the subgroup of $B$ fixing $R$ pointwise (in this notation, a line is also considered to be a point set).\\

Generalized quadrangles and polygons were introduced by Tits in a famous work on triality \cite{Ti},  in order to propose an axiomatic and combinatorial treatment for 
semisimple algebraic groups (including Chevalley groups  and groups of Lie type) of relative rank 2. The most important  generalized polygons are those with gonality $4$ | the {\em generalized quadrangles}. Standard reference works are \cite{PT}, \cite{SFGQ, LEGQ} and \cite{TGQ}. Note that projective planes are nothing else than generalized $3$-gons. Also, the class of generalized polygons coincides with the class of Tits buildings of rank $2$; so they are the rank $2$ residues for Tits buildings of higher rank. \\

Let us proceed with the model theoretic approach to locally finite quadrangles, first initiated by Gregory Cherlin in \cite{Cher}.

\bigskip
\section{Indiscernibles and generalized quadrangles}
\label{indis}

Let $\Upgamma$ be a generalized quadrangle, and let $(S,\leq)$ be a totally ordered set.  A set $\mL = \{ M_s\ \vert\ s \in S \}$ of lines is {\em indiscernible} over $(S,\leq)$ if for any
two increasing sequences $M_{s_1},M_{s_2},\ldots,M_{s_n}$ and $M_{s_1'},M_{s_2'},\ldots,M_{s_n'}$ (of the same length $n$) of lines in $\mL$, there 
is an automorphism $\upgamma$ of $\Upgamma$ mapping $M_{s_i}$ onto $M_{s_i'}$ for each $i$. Here, by ``increasing'' we mean that in $S$, $s_1 < s_2 < \ldots < s_n$ and $s_1' < s_2' < \ldots < s_n'$. The type of indiscernibility used here is also called {\em automorphism-theoretic indiscernibility}. Note that a priori, we do not ask that $\upgamma$ globally fixes $\mL$. \\

% It is indiscernible {\em over $D$}, if $D$ is a finite set of points and lines fixed by the automorphisms just described.\\

By combining the Compactness Theorem and Ramsey's Theorem \cite{Hodge} (in a theory which has a model in which a given definable set is infinite), one can prove the following.

\bt[G. Cherlin \cite{Cher}]
\label{Cherlem}
Suppose there is an infinite locally finite generalized quadrangle with finite lines. Then there is an
infinite locally finite generalized quadrangle $\Upgamma$ containing an indiscernible sequence $\mathcal{L}$ of parallel (= mutually skew) lines, of any 
specified order type. 
%The sequence may be taken to be indiscernible over the set $D$ of all points incident with one fixed line $L$ of $\Upgamma$.
\eop 
\et

%Clearly, $D$ may supposed to be fixed pointwise in Theorem \ref{Cherlem}.

\br
\label{pol}
{\rm
%\begin{itemize}
%\item[{\rm(i)}]
%Cherlin states the theorem only for generalized quadrangles, but the same statement also holds for ``general'' generalized quadrangles.
%\item[{\rm(ii)}]
In Theorem \ref{Cherlem}, $\Upgamma$ may supposed to be generated by $\mathcal{L} \cup \{L\}$ \cite{Cherpriv} (where $L$ is seen as a point set). ``Generated by $\mathcal{L} \cup \{L\}$'' here, means 
that $\Upgamma$ is the only subquadrangle which contains $L$, the lines of $\mL$, and the points of $L$. }
%\end{itemize}}
\er

%In this paper we will show that infinite locally finite generalized $4$-gons do not exist. By Theorem \ref{Cherlem},  it will be sufficient to
%prove the result for $t$ countable. (We will make this claim more precise in the next section.)

\subsection*{Moorhouse's text}

We strongly advise the reader who is interested in understanding the model theoretic details behind Cherlin's paper \cite{Cher}, to consult the excellent 
notes of G. Eric Moorhouse \cite{Moor}. Moorhouse's text contains all the details a geometer needs to start on this topic.

\bigskip
\section{Some further conventions}

Let $\Upgamma$ be a locally finite GQ, and let $X$ be a set of points or lines (not both), which is indiscernible over the totally ordered set $(S,\leq)$. By $A$, we denote 
the subgroup of $\Aut(\Upgamma)$ which globally stabilizes $X$.  So $A = \Aut(\Upgamma)_X$.

Let $(S,\leq)$ be a totally ordered set, and let $S'$ be a subset of $S$; then by $(S',\leq)$ we denote the totally ordered set which arises by restricting $\leq$ to $S'$.

%If $S \subseteq \mL$, by $\Upgamma(S)$ we denote the full subquadrangle generated by $S \cup \{L\}$ in the sense of Remark \ref{pol}(ii) (since $L$
%is fixed throughout, we do not specify $L$ in the notation). So $\Upgamma = \Upgamma(\mL)$.\\

%Call a sequence of lines $N_{i_1},N_{i_2},\ldots,N_{i_l}$ in $\mL$ {\em increasing} if $i_l < i_{l'}$ for 
%$l < l'$.

%\bl
%\label{indisc}
%If $S \subset \mL$ is a subset of $\mL$, then $\Upgamma(S)$ is disjoint from any line of 
%$\mL \setminus S$.
%\el
%{\em Proof}.\quad
%Immediate from indiscernibility.
%\eop \\

 %Let $\mathrm{Sub}(\Upgamma)$ be the set of subGPs of $\Upgamma$.

%\bl
%The map 
%\[ \Psi: 2^{\mL} \mapsto \mathrm{Sub}(\Upgamma): S \mapsto \varphi(S) = \Upgamma(S)               \]
%is an injection. 
%\el
%{\em Proof}.\quad
%Immediate.
%\eop \\

%\medskip
%The following lemma is folklore, and easy to prove. (In its statement, ``countable'' also comprises
%the finite case.)

%\bl
%\label{lemcount}
%Let $\Delta$ be a generalized polygon, and $K$ be a subset of points of countable size.
%Let $\Delta(K)$ be the subpolygon of $\Delta$ generated by $K$. Then the number of points and 
%lines of $\Delta(K)$ is countable.
%\eop \\
%\el

%\bc
%We have that $t$ is countable, as is the number of points and lines of $\Upgamma$. 
%\ec

%{\em Proof}.\quad
%Since $\Upgamma$ is generated by $\{L\} \cup \mL$, it is generated by a countable number of points.
%By Lemma \ref{lemcount}, the number of points and lines is countable. Clearly $t$ also is.
%\eop \\

\bigskip
\section{Moorhouse's questions}
\label{Moorques}

Let $\mS$ be a GQ, and let $X$ be a line set in $\mS$. If no two distinct elements in $X$ intersect, then we call $X$ a {\em partial spread}. If every point of $\mS$ 
is incident with a line of $X$, then we call $X$ a {\em spread}. Dually, we speak of {\em partial ovoids} and {\em ovoids}.

\begin{quote}
\begin{mdframed}[linewidth=2.5,linecolor=orange, topline=false,rightline=false,bottomline=false]
\item[(IND$_{ov}$)]
Suppose $\Upgamma$ is a locally finite GQ of order $(\upomega,k)$ with $k$ finite, which contains an ovoid $\mO$, which is indiscernible over some totally ordered 
set $(S,\leq)$. What can we say about $\Upgamma$? 
\item[(IND$_{sp}$)]
Suppose $\Upgamma$ is a locally finite GQ of order $(\upomega,k)$ with $k$ finite, which contains a spread $\mL$, which is indiscernible over some totally ordered 
set $(S,\leq)$. What can we say about $\Upgamma$? 
\end{mdframed}
\end{quote}

Relative to the assumption that $\Upgamma$ has order $(\upomega,k)$, we will see that (IND$_{ov}$) is easier than (IND$_{sp}$). Both questions will 
have the same answer, though: \ul{such quadrangles $\Upgamma$ cannot exist}. \\

%In this paper, we consider all generalized quadrangles $\Upgamma$ for which the totally ordered set $(S,\leq)$ is countable, that is, for which $\mO$ (case IND$_{ov}$) and $\mL$ (case IND$_{sp}$) have a countable number of elements. 

\bigskip
\section{Quadrangles with an indiscernible ovoid or spread}
\label{HOW}

Let $\Upgamma$ be a thick generalized quadrangle of order $(\upomega,k)$, with $k$ finite and $\upomega$ not finite. %So if we want to apply (EMB) in the setting of hypothesis (IND$_{sp}$), $\overline{\Upgamma}$ is the unique  ...

\subsection{Hypothesis (IND$_{ov}$)}

First suppose that $\Upgamma$ satisfies (IND$_{ov}$): we let $\mO$ be an ovoid of indiscernible points, over the totally ordered set $(S,\leq)$. 
Let  $w$ be a point which is not in $\mO$, and let $X, Y, Z$ be three different lines incident with $w$; since $\mO$ is an ovoid, each of these 
lines intersects $\mO$ in precisely one point, say (respectively) $x, y$ and $z$. Without loss of generality, we suppose that 
\begin{equation}
x < y < z.
\end{equation}   

For each $\upepsilon \in \mO$ for which $\upepsilon > z$, there is an automorphism of $\Upgamma$ which fixes $x$ and $y$, and sends $z$ to $\upepsilon$. As it follows that each of the points of
the set $\{ \upepsilon\ \vert\ \upepsilon > z\ \mbox{and}\ \upepsilon \in \mO \}$ is collinear with some point of $\{ x,y \}^{\perp}$, the fact that $\vert \{x,y\}^{\perp} \vert = k + 1 \in \mathbb{N}$ leads to the fact that 
some point of $\{ x,y \}^{\perp}$ must be incident with an infinite number of lines, contradiction. 

This solves the first question. \eop \\

Note that the proof does not use the fact that $S$ is countable. So in general, case (IND$_{ov}$) is resolved (through the most general notion of indiscernibility available).

\subsection{Hypothesis (IND$_{sp}$)}

We now turn to the second hypothesis, Hypothesis (IND$_{sp}$). This case is more delicate, and the previous proof does not apply. In fact, we will need a stronger version of indiscernibility, which we will describe below. 

{ 
\subsection*{(Strong) Automorphism-theoretic indiscernibility}

From this point on, if we use the term 
``indiscernible," we will mean {\em strong} automorphism-theoretic indiscernible. This means that the automorphisms of $\Upgamma$ required to make $\mL$ indiscernible over $(S,\leq)$, \ul{stabilize the set $\mL$}. 
We will also assume that all automorphisms of $(S,\leq)$ are induced by elements of $\Aut(\Upgamma)_S$; we do not need the full force of this assumption, but it is a natural property from the incidence-geometrical point of view nevertheless.   
Although many of our arguments work for the more general notion of automorphism-theoretic indiscernibility, unfortunately not all of them might do.} \\
%In fact, the only thing we need is that there exist two different elements $a < b$ in $S$ and an automorphism $\upeta$ of $(S.\leq)$ which does not fix all points of $S$ in $(a,b)$, which fixes every point of $S$ not in $(a,b)$, and which is induced (in the obvious way) by some automorphism $\widehat{\upeta}$ of $\Upgamma$. Such automorphisms $\upeta$ obviously exist if we ask that {\em each} automorphism of $(S,\leq)$ is induced.} \\

Let $\mL$ be a spread of lines of $\Upgamma$, indiscernible over the totally ordered set $(S,\leq)$. We do not restrict ourselves to the countable case. A first direct observation is that directly due to (IND$_{sp}$), we know that $(S,\leq)$ is a dense linear order without end points (so that by a result of Cantor, 
$(S,\leq) \simeq (\mathbb{Q},\leq)$ in the countable case). 

Let $u$ be any point in $\Upgamma$, and let $V$ be any line incident with $u$, not contained in $\mL$. Let $S_V \subset S$ be the set of elements $s$ in $S$ 
which correspond to the lines in $\mL$ that meet $V$.  

Now introduce an equivalence relation on $S_V$ as follows: $a \sim b$ if only a finite number (including $0$) of elements in $S_V$ are contained between $a$ and $b$. Denote the set of equivalence classes in $S_V$ by $\widetilde{S_V}$, and note that each equivalence class trivially has a finite or countable number of elements. Also, $\leq$ naturally induces an order on $\widetilde{S_V}$, and we will keep using the same symbol ``$\leq$'' for that order. Putting $\widehat{S_V} := S \setminus S_V$, we similarly introduce $\widetilde{\widehat{S_V}}$. 

Now note the following property:
\begin{quote}
$(\widetilde{S_V},\leq)$ is a dense total order. 
\end{quote}
(Indeed, if $r, s \in \widetilde{S_V}$ and $r < s$, then if no $t \in \widetilde{S_V}$ would exist for which $r < t < s$, we would have that $r = s$.)\\

We will call elements of $\widetilde{S_V}$, resp. $\widetilde{\widehat{S_V}}$,  ``singletons'' if they only contain one element as an equivalence class in $S_V$, resp. $\widehat{S_V}$. We distinguish a number of possibilities, depending on properties of ``non-singletons.''\\

\subsection*{How to proceed}

Suppose $\overline{S_V}$ is a subset in $S$ such that there is some element $\upalpha \in A$ for which $\upalpha(S_V) = \overline{S_V}$, and such that 
$S_V \subset \overline{S_V}$ (strict inclusion). Then we have obtained a contradiction, since $\{ V,V^\upalpha \}^{\perp} = \{ L_s \ \vert\ s \in S_V \}$, while $V^\upalpha$ is incident with points which are not incident with lines of 
$\{ V, V^\upalpha \}^{\perp}$. { (For, as $S_V \subset \upalpha(S_V)$, $V^\upalpha$ meets all $L_t$ for $t \in S_V$. Take $s \in S_V$ with $s' = s^\upalpha \not\in S_V$. Then $V^\upalpha$ meets $L_{s'}$ in a point $x$. There is a point $y$ on $V$ collinear with $x$. It follows that $y \I L_r$ for some $r \in S_V$, so $L_r$ is a second line through $y$ meeting $V^\upalpha$, contradiction.)} 

In each of the next paragraphs, we will search for such automorphisms $\upalpha$. Only if $S_V = S$ for all choices of $u$ and $V$, no contradiction arises. 
And this is precisely the case where $\Upgamma$ is a grid.\\

Below, we say that $U$ is a {\em dense subset} of the linear order $(V,\leq)$ if 
for all $a < b$ in $V$, there is a $c \in U$ such that $a < c < b$. \\

{ In what follows, we assume without loss of generality that there are no nontrivial intervals contained in $\widehat{S_V}$; if intervals would occur both with only  points in $S_V$ and only points in $\widehat{S_V}$, we proceed as follows. 
Let $(a,a')$ be an open interval for which all interior elements are in $\widehat{S_V}$, 
and let $(q,q')$ be an open interval for which all the interior elements are in $S_V$. 
We suppose without loss of generality that $q < q' < a < a'$ (the other case is handled by a symmetrical argument). 

Now define an automorphism $\upalpha \in A$ with the following properties: 

\begin{itemize}
\item
$\upalpha(q) = q$ and $\upalpha(a') = a'$;
\item
if $r < q$ or $r > a'$, then $\upalpha(r) = r$ ($r \in S$);
$\upalpha(q') \in (a,a')$;
\item
$\upalpha(q') < \upalpha(a) < a'$.   
\end{itemize}

It is easy to see that such $\upalpha$ exist in $A$, and it is obvious that $\upalpha(S_V)$ strictly contains $S_V$.
Now proceed as in the beginning of this section. \\

Also, if additionally there would be 
no intervals with only interior points in $S_V$, then both $S_V$ and $\widehat{S_V}$ are dense subsets of $S$, 
and we can proceed with a particular case of the Back \& Forth argument at the end of the next section.} \\

\medskip
\section*{There are only singletons (in both $\widetilde{S_V}$ and $\widetilde{\widehat{S_V}}$)}

This means that $S_V = \widetilde{S_V}$ (modulo bracket notation), and that $(S_V,\leq)$ is a dense linear order, possibly with boundaries. Moreover, $(\widehat{S_V},\leq)$ 
is also a dense linear order (with $\widehat{S_V} := S \setminus S_V$). 

%Now define a new linearly ordered set $(\widetilde{S},\leq)$, by introducing the following equivalence relation: $x \sim y$ if either $x, y \in S_V$ and 
%both are contained in the same $S_V$-interval, or $x, y \in \widehat{S_V}$ and both are contained in the same $\widehat{A}$-interval. 

Let $q < q'$ be points in $S$; then we write $q \asymp q'$ if each point $r \in S$ for which $q < r < q'$, is contained in $S_V$. If $q = q'$, then also $q \asymp q'$. 
We will assume without loss of generality that ``$\asymp$'' defines 
an equivalence relation on $S$ (and $S_V$), and we denote the equivalence class containing some element $v \in S$ by $[v]$. For $\asymp$ not to be 
an equivalence relation, we need a violation of transitivity, and hence we must have distinct points $a, b, c$ in $S$ such that $a \asymp b$ and $b \asymp c$ while $a \not\asymp c$. 
But that case can be handled as the case in section \ref{box3}. 
If $(u,v)$ is an open 
interval in $S$, then we will, at various points, also consider the induced equivalence classes of $\asymp$ inside $(u,v)$. 
Call points in $S$ (or in $(u,v) \subseteq S$) {\em isolated} if their class only has one element, and of ``type 2'' otherwise. In the latter case, an equivalence class 
$[v]$ is an interval of which the endpoints cannot be both elements of $\widehat{S_V}$ (by the main assumption in this section). 
We have a number of different types of such intervals: 
\begin{itemize}
\item
closed with both endpoints in $S_V$; 
\item
half-open $(b,a]$ or $[a,b)$ with $a \in S_V$ and $b \in \widehat{S_V}$; 
\item
closed with one endpoint in $S_V$ and the other in $\widehat{S_V}$. 
\end{itemize}
Each of these is considered as a different type 2 equivalence class (but we will still say that their ``main type'' is 2). Usually we do not  consider intervals with non-finite endpoints (they are usually not needed in our arguments). The type of an interval which corresponds to an equivalence class, is by definition the type of its class. For the sake of convenience, denote the set of all types with main type $2$ by $\mT$. 

Let $a < c$ points in $S$. By the expression ``$a < \bI < c$,''  
with $\bI$ an interval of type $\bT \in \mT$, we mean that $a < e$ and $e' < c$, with $e < e'$ the endpoints of $\bI$. \\

%Define the following property, with $\mT' \subseteq \mT$ and $S' \subseteq S$:
%\begin{quote}[DEN($\mT'$,$S'$)]
%for all $q_1 < q_2$ with $q_1, q_2$ elements in $S$, and for each type $\bT \in \mT'$, there is an element $r$ of type $\bT$, such that $q_1 < r < q_2$. Also, for each $q_1 \in S$, there is such an element $r$ for which $r < q_1$, and 
%for each $q_2 \in S$, there is such an element $r$ for which $q_2 < r$. 
%\end{quote}

%We will only consider the first property; the others have similar proofs. Below, we show that ...
%We consider the possible types. 

We proceed as follows. 

%First consider some interval $(q,q')$, and let $\bI_1$ and $\bI_2$ be two distinct intervals of main type 2 in $(q,q')$; this means they arise from $\asymp$ restricted to $(q,q')$.  Suppose they are maximal in $(q,q')$ with respect to being of main type 2. Suppose first that they share an endpoint $e_2 = e_1'$, where 
%$e_2$ is the largest endpoint of $\bI_1$ and $e_1'$ is the smallest endpoint of $\bI_2$. Then since both intervals are maximal, we have that $e_2$ cannot be in $S_V$; so it is a point of $\widehat{S_V}$. In that case, proceed as in subsection \ref{box3} (where $e_2$ is the point that is moved).  

%Henceforth, we may and will suppose that if such distinct intervals 
%exist, they do not share endpoints. \\

\quad(A)\quad
If there exist $q_1 < q_2$ in $S$ such that $q_1 \not\asymp q_2$, and 
such that for no element $\bI$ of type $\bT$ (with main type 2),  we have that $q_1 < \bI < q_2$, then we replace $S$ by $(q_1,q_2)$. Now repeat the same argument with $S$ replaced by $(q_1,q_2)$, and $\bI$ replaced by an element with type in $\mT \setminus \{ \bT\}$. 
After having considered all possible types,  we end up with an open interval $(q_1^*,q_2^*)$ (which might still equal $S$ but which otherwise has a finite endpoint) \ul{with $q_1^* \not\asymp q_2^*$, and 
such that for all $q_1 < q_2$ in this interval which are not contained in $[q_1^*] \cup [q_2^*]$ and for which $q_1 \not\asymp q_2$, we have an interval $\bI$ of type $\bT \in \mT'$ such that $q_1 < \bI < q_2$.} Here $\mT'$ is a subset of $\mT$ of types with main type 2; in principle, the set $\mT'$ can be empty. In that case, there are no intervals of type $\bT \in \mT$ inside 
$(q_1^*,q_2^*)$, except, possibly, intervals containing $q_1^*$ or $q_2^*$. This case is handled as follows:
\begin{itemize}
\item
if there is exactly one point of $\upomega \in \widehat{S_V}$ contained in $(q_1^*,q_2^*)$, then proceed as in section \ref{box3}.
\item
if we are not in the previous case, then $\Upomega := (q_1^*,q_2^*) \setminus [q_1^*] \cup [q_2^*]$  consists only of isolated points of $S_V$ and 
$\widehat{S_V}$, and both sets are dense in $\Upomega$. The same Back \& Forth argument which will be carried out at the end of this section 
then ends the proof.   
\end{itemize}

%If no $r$ exists in $\widetilde{S_V}$ of type 2, then in $(q_1,q_2)$, $S_V \cap (q_1,q_2)$ and $\widehat{S_V} \cap (q_1,q_2)$ are dense subsets. Now we can easily construct the desired $\upalpha$, since $\Aut(\mathbb{Q},\leq)$ acts transitively on the pairs $(U,V)$ with $U$ a dense subset of $\mathbb{Q}$, and $V = \mathbb{Q} \setminus U$ a dense subset in $\mathbb{Q}$ (apply this property on $(q_1,q_2)$, and let $\upalpha$ fix all elements $s \in S$ outside this interval).

\quad(B)\quad
Now suppose there is no isolated point of $S_V$  in some $(q_1,q_2) \subseteq (q_1^*,q_2^*)$ with $q_1 \not\asymp q_2$. Note that  intervals of type 2 must exist inside $(q_1,q_2)$, since we otherwise have nontrivial intervals whose interior points are solely in $\widehat{S_V}$. Replace $(q_1^*,q_2^*)$ by $(q_1,q_2)$ (while keeping the same notation). If there would be $q_1' < q_2'$ inside $(q_1^*,q_2^*)$ with $q_1' \not\asymp q_2'$ and no isolated point of $\widehat{S_V}$ in between, then we replace once again 
$(q_1^*,q_2^*)$ by $(q_1',q_2')$ while keeping the same notation; otherwise, the isolated points of $\widehat{S_V}$ define a dense set inside the equivalence classes of 
$\asymp$ in $(q_1^*,q_2^*)$. Both cases can be handled in the same way, as we will see below. 
Let $\bI = [u,v]$ or $[u,v)$ be an equivalence class in $(q_1^*,q_2^*)$ such that it does not contain $q_1^*$; without loss of generality, we 
will suppose that $u \in S_V$ (all the other variations are completely similar as this case). Now let $\upalpha$ be an order automorphism of $S$ such that for all $r \geq v$ and $r \leq q_1^*$, we have that,  $\upalpha(r) = r$, and such that  $u \not\asymp \upalpha(u) < u$ be such that $\upalpha(u) \in (q_1^*,q_2^*)$. 
%Define 
%\begin{equation}
%\mD\ :=\ \{ \upalpha'(t)\ \vert\ t \in \widehat{S_V}, t \in (q_1^*, u) \}.
%\end{equation}
%Enumerate the elements of $\mD$ as $d_1, d_2, \ldots$. 
We also require that 
\begin{equation} 
\upalpha:\ (q_1^*,u)\ \mapsto\ (q_1^*,\upalpha(u))
\end{equation}
sends $\widehat{S_V} \cap (q_1^*,u)$ surjectively to $\widehat{S_V} \cap (q_1^*,\upalpha(u))$, such that for $t \in \widehat{S_V} \cap (q_1^*,u)$ we have that $\upalpha(t)$ is a (left resp. right) boundary point if $t$ is a (left resp. right) boundary point | by our definition of $(q_1^*,q_2^*)$, we can easily construct $\upalpha$ through a standard Back \& Forth argument. (A similar detailed Back \& Forth argument will be carried out at the end of this section.)

%We will construct a new order automorphism $\upalpha$ by first defining it on the points  $t \in \widehat{S_V} \cap (q_1^*,u)$. It is sufficient to define the action $d \mapsto \widetilde{d}$ on the points of $\mD$. Suppose $m$ is the first index for which $d_m = \upalpha'(t)$ is not contained in $\widehat{S_V}$. Suppose 
%$\widetilde{d_a} < d_m < \widetilde{d_b}$ for the images $\widetilde{d_a}$ and $\widetilde{d_b}$ which are both contained in $\widehat{S_V}$. Take an element $s_m \in \widehat{S_V} \cap (\widetilde{d_a},\widetilde{d_b})$ which is isolated if $t$ is isolated, and a (left resp. right) boundary point if $t$ is a (left resp. right) boundary point | by our definition of $(q_1^*,q_2^*)$, this is always possible. Now put $\upalpha(t) = \widetilde{d} =  s_m$. Put $\upalpha = \upalpha'$ for all 
%$r \leq q_1^*$ and $r \geq u$, and put $\upalpha(t) = s_m$ for all $t \in \widehat{S_V} \cap (q_1^*,u)$ with $d_m = \upalpha'(t)$. One now completes the construction of $\upalpha$ by defining images for  the points in $S_V \cap (q_1^*,u)$ while respecting the order of the already defined images, and this is done through a simple Back \& Forth argument. 

For such $\upalpha$, we have that $\upalpha(S_V) \supset S_V \ne \upalpha(S_V)$.   

So in what follows, we may and will \ul{assume that for all $q_1 < q_2$ in $(q_1^*,q_2^*)$ with $q_1 \not\asymp q_2$, there is an $r \in S_V$ with $[r] = \{r\}$ such that $q_1 < r < q_2$.} \\

\quad(C)\quad
Suppose that $q_1 < q_2$ in $(q_1^*,q_2^*)$ are such that $q_1 \not\asymp q_2$.  
Then we may assume (through an argument similar as in the previous point) that we have points $r$ in $\widehat{S_V}$ for which $\{ r\} = [r]$ such that $q_1 < r < q_2$. \\
% \item
%If  $q_1 < q_2$ are elements in $(q_1^*,q_2^*)$ such that $q_1 \not\asymp q_2$, and there is no point $r$ of ${S_V}$ in between, then we immediate 
%have a contradiction as then $q_1 \asymp q_2$.  

Note that for all $q_2 \in (q_1^*,q_2^*)$ for which 
$q_2 \not\asymp q_2^*$, we can find isolated points $r$ in $S_V$, resp. $\widehat{S_V}$, resp. intervals $\bI$ of type $\bT$ in $\mT'$ such that 
$q_2 < r < q_2^*$, resp. $q_2 < r < q_2^*$, resp. $q_2 < \bI < q_2^*$. Also, we have that for all $q_1 \in (q_1^*,q_2^*)$ for which 
$q_1^* \not\asymp q_1$, we can find isolated points $r$ in $S_V$, resp. $\widehat{S_V}$, resp. intervals $\bI$ of type $\bT$ in $\mT'$ such that 
$q_1^* < r < q_1$, resp. $q_1^* < r < q_1$, resp. $q_1^* < \bI < q_1$. Note that if $\bI_1 < \bI_2$ ($\bI_1, \bI_2$ intervals), then for $q_1 \in \bI_1$ and $q_2 \in \bI_2$, we have that $q_1 \not\asymp q_2$ and the above properties apply.

\medskip
We now construct the desired $\upalpha \in A$.  \\

{\bf BACK \& FORTH argument}.\quad
{\em First suppose the number of equivalence classes of $\asymp$ is countable.}\quad 
Enumerate the equivalence classes of $\asymp$ in $(q_1^*,q_2^*)$: $s_1,s_2,\ldots$. We do not include the classes of $q_1^*$ and/or $q_2^*$ in the enumeration; below, once we define $\upalpha$, we put $\upalpha(L_r) = L_r$ for each element $r$ 
inside $[q_1^*]$ and/or $[q_2^*]$. 

Now let $r$ be  an element in $(q_1^*,q_2^*) \cap \widehat{S_V}$ with $r$ not in $[q_1^*]$ nor $[q_2^*]$ and such that $[r] = \{r\}$, and define: $S_V' := S_V \cup \{ r\}$, 
$\widehat{S_V'} := \widehat{S_V} \setminus  \{ r\}$. Relative to these new sets $S_V'$ and $\widehat{S_V'}$, we define $\asymp'$ similarly as $\asymp$, and adapt the definition of types. Note that $[r]' = \{r \}$. 
Enumerate the equivalence classes of $\asymp'$ as $s_1',s_2',\ldots$.\\

Let $(\mathcal{E},\leq)$ be the set of equivalence classes of $\asymp'$ in $(q_1^*,q_2^*)$ endowed with the induced order relation; the reason why the Back \& Forth argument below to construct $\upalpha$ works, is the fact that the equivalence classes of each fixed type (including singletons) are dense in $\mathcal{E}$. \\

Suppose that $m$ is the smallest index such that $s_m$ is not yet paired with any member of $s_1',s_2',\ldots$. Let $s_k'$ be an element that is not yet paired, and 
which is \ul{of the same type as $s_m$}, while respecting the order of the elements. Pair $s_m$ to $s_k'$. 
Now let $n$ be the smallest index for which $s_n'$ is not yet paired, and pair it similarly 
with an element of $s_1,s_2,\ldots$ that is not yet paired. One uses the density properties of the previous paragraph to show that the pairings can be done in each step without violating the order. 
$\Big($For instance, if $s_{l_1} < s_m < s_{l_2}$ with $l_1, l_2$ indexes which were already taken care of in previous steps, and $s_m = \{\rho\}$ with $\rho$ in $\widehat{S_V}$, 
and $s_{l_1}$ was paired to $s_{r_1}'$ and $s_{l_2}$ to $s_{r_2}'$, then $s_{r_1}' < s_{r_2}'$ and we can find a singleton $s_n'$ in $\widehat{S_V}$ such that 
$s_{r_1}' <  s_n' < s_{r_2}'$ ($s_n' = \{ \rho'\}$, $\rho' \in \widehat{S_V}$). That is the point to which $s_m$ is paired. All the other cases are similar.$\Big)$ \\

As such, we have constructed a bijection between $s_1,s_2,\ldots$ and $s_1',s_2',\ldots$; call it $\underline{\upalpha}$. Now we define the desired element of $A$ as follows: 
$\upalpha$ acts trivially on lines indexed by elements outside $(q_1^*,q_2^*)$; if $[s_i] = \{s_i\}$, then $\upalpha(L_{s_i}) := L_{\underline{\upalpha}(s_i)}$; if $s_j$ is an interval, then choose an arbitrary order bijection $\widetilde{\upalpha}$ 
between $s_j$ and $\underline{\upalpha}(s_j)$ which also defines  the images for the endpoints 
(such bijections exist since $s_j$ and $\underline{\upalpha}(s_j)$ are of the same type), and set $\upalpha(L_s) = L_{\widetilde{\upalpha}(s)}$ for the points $s$ in $s_j$, and for its endpoints.  \\

Then $\upalpha \in A$ is such that $\upalpha(S_V) = S_V' = S_V \cup \{r \}$.  \\

{\em Now suppose that the number of equivalence classes of $\asymp$ is not countable.}\quad In that case, we apply transfinite recursion to construct $\underline{\upalpha}$.
We index the equivalence classes of $\asymp$ and $\asymp'$ by ordinals in the usual way (and note that both sets have the same size). Suppose $s_\upalpha$ is an element which was already paired to $s'_\upalpha$; then below, $\widehat{s'_\upalpha}$ denotes the 
sequence $\Big\{ s'_\epsilon \Big\}_{\epsilon \leq \upalpha}$. 
For equivalence classes of $\asymp$ which are indexed by successor ordinals and not yet paired, we refer to the proof of the countable case. If the index $\upalpha$ is a limit ordinal, then we define 
\begin{equation}
\widehat{s'_\upalpha} \:=\ \bigcup_{\upbeta < \upalpha}\widehat{s'_\upbeta}. 
\end{equation}
As such, we have defined $\underline{\upalpha}$. The rest is the same as in the countable case. \\

%Enumerate the equivalence classes of $\asymp$ in $(q_1^*,q_2^*)$: $\{ s_\upalpha\ \vert\ \upalpha < \upbeta \}$, with $\upbeta$ the smallest ordinal for which $\vert \Epsilon \vert = \upbeta$. (Just like in the countable case, we do not include the classes of $q_1^*$ and/or $q_2^*$ in the enumeration; and again, once we define $\upalpha$, we put $\upalpha(L_r) = L_r$ for each element $r$ inside $[q_1^*]$ and/or $[q_2^*]$.)  

{\em Alternative argument} (in both the countable and uncountable case).\quad 
One can also define $\upalpha$ through the method presented in the penultimate item (B) in the enumeration above (which essentially solely concentrates on the images 
of $S_V$, starting from an order automorphism $\upalpha'$ as above). The essence is the same. \\

%It is now easy to construct such $\upalpha$ in $A$, as $A$ acts transitively on the dense subsets in $\mathbb{Q}$ with dense complements. 

\medskip
\section*{There is at least one non-singleton { (in $\widetilde{S_V}$ or $\widetilde{\widehat{S_V}}$)} which is not finite}

Call it $s$. { Note that for each element $r$ in this class, there exists are elements $a, b \in S$ with $a < r < b$, such that 
$(a,b)_S \cap s = \{r \}$.}  

Now we construct an element of $A$. 

First of all, let $\upgamma(t) = t$ for all $t \in S$ such that $t < s$ (by which we mean that $t < s'$ for all $s' \in s$). Similarly, let 
$\upgamma(t) = t$ for all $t > s$. Now consider consecutive elements $u, v$ in $s$, and let $w \in S$ be strictly contained 
between $u$ and $v$. Define $s' = s \cup \{w\}$. Then $(s,\leq) \simeq (s',\leq)$, and we assume that $\upgamma$ induces 
an order isomorphism between $(s,\leq)$ and $(s',\leq)$. 

Now we extend $\upgamma$ to an element $\overline{\upgamma} = \upalpha \in A$ (by a standard Back \& Forth argument). Then $\overline{S_V} := \overline{\upgamma}(S_V)$ 
strictly contains $S_V$. Now proceed as before. \\

\begin{remark}{\rm 
The reasoning which we will see below in subsection \ref{box3} can also be applied to the case considered  in this section. 
}
\end{remark}

\section*{There is at least one non-singleton { (in $\widetilde{S_V}$ or $\widetilde{\widehat{S_V}}$)}, and all of them are finite}

We consider two sub-cases. 

\subsection{At least one non-singleton has size $\geq 3$}
\label{box3}

Let $s$ be such a non-singleton { which we assume to be in $\widetilde{S_V}$ without loss of generality}, and let $r$ be a point of $s$ which is not the largest nor smallest element of $s$. Now consider an element $\upalpha$ in $A$ 
such that all elements of $S_V \setminus \{r\}$ are fixed, and such that $r$ is moved inside $s$. Such $\upalpha$ are easily seen to exist, since $r$ has a positive distance to each of its neighbors in $S_V$. 

Then $V^\upalpha$ necessarily is a line different from $V$, and $\{ V, V^\upalpha \}^{\perp}$ contains all the lines of $\mL$ which are indexed by an element of $S_V \setminus \{r \} = \upalpha(S_V) \setminus \{\upalpha(r)\}$. The same property is true for all elements in the infinite set 
\begin{equation}
\upsigma(V) := \{ \upalpha^m(V)\ \vert\ m \in \mathbb{N} \setminus \{0\} \}.
\end{equation}

{It is straightforward to see that there is one extra line in $\Upgamma$ which meets all lines of $\{ V \} \cup \upsigma(V)$,  containing the points 
$\widetilde{r}$ of $V$, where $\{ \widetilde{r} \} = V \cap L_r$ and $L_r$ is the line of $\mL$ indexed by $r$,  $\upalpha(\widetilde{r})$ of $\upalpha(V)$, $\upalpha^2(\widetilde{r})$ of $\upalpha^2(V), \ldots$ (because otherwise triangles arise).
Hence 
\begin{equation}
\{V \} \cup \upsigma(V) \subseteq \{ V, V^{\upalpha} \}^{\perp\perp}. 
\end{equation}
}

Now consider a point $\upomega$ not contained in the set of points incident with lines in $\{ V \} \cup \upsigma(V)$; by projecting $\upomega$ on the lines in the latter line set, we immediately obtain a contradiction since there are only a finite number of lines incident with $\upomega$. (If $U, W$ are different lines in $\{V \} \cup \upsigma(V)$, then the projection line from $\upomega$ to $U$ is different than the projection line from $\upomega$ to $W$, { so since $\{ V \} \cup \upsigma(V)$ is infinite, so must be the set of lines incident with $\upomega$.})

\subsection{All non-singletons have size $2$}

We need to consider two sub-cases.

\subsubsection{Both $\widetilde{S_V}$ and $\widetilde{\widehat{S_V}}$ have non-singletons}
\label{boxes}

{
Let $(a,a')$ be an open interval for which all interior elements are in $\widehat{S_V}$, 
and let $(q,q')$ be an open interval for which all the interior elements are in $S_V$. By assumption, these 
intervals exist. This case was already resolved in section \ref{HOW} (in the paragraph ``How to proceed''). }

%We suppose without loss of generality that $q < q' < a < a'$ (the other case is handled by a symmetrical argument). 

%Now define an automorphism $\upalpha \in A$ with the following properties: 

%\begin{itemize}
%\item
%$\upalpha(q) = q$ and $\upalpha(a') = a'$;
%\item
%if $r < q$ or $r > a'$, then $\upalpha(r) = r$ ($r \in \mathbb{Q}$);
%$\upalpha(q') \in (a,a')$;
%\item
%$\upalpha(q') < \upalpha(a) < a'$.   
%\end{itemize}

%It is easy to see that such $\upalpha$ exist in $A$, and it is obvious that $\upalpha(S_V)$ strictly contains $S_V$. 
%Now proceed as before. Note that if $q = -\infty$ and/or $a' = +\infty$, the approach above trivially adapts just fine.} 

\subsubsection{Either $\widetilde{\widehat{S_V}}$ or $\widetilde{S_V}$ has no non-singletons}

This case can be handled through a similar Back \& Forth method as the case without non-singletons. \\

%{ 
%\subsubsection{$S_V$ has no non-singletons}

%This case is similarly handled as in the previous case.}\\

The finally completes the proof of the answers to both questions. \eop \\

%\medskip
%{\bf FORTHCOMING}.\quad 
%In a forthcoming sequel to this paper \cite{part2} (which we hope to have finished rather sooner than later), we will handle the uncountable version of (IND$_{sp}$). 

\newpage
%\vspace*{1cm}

\end{document}